\documentclass[12pt]{amsart}
\usepackage[centertags]{amsmath}
\usepackage{amsfonts}
\usepackage{amssymb}
\usepackage{amsthm}
\usepackage{newlfont}
\usepackage{bm}
\newlength{\defbaselineskip}
\newcommand{\setlinespacing}[1]%
{\setlength{\baselineskip}{#1 \defbaselineskip}}

\def\CC{{\mathcal C}}
\def\DD{{\mathcal D}}

\def\MM{{\mathcal M}}

\def\RR{{\mathcal R}}

\def\TT{{\mathcal T}}

\def\bbC{\mathbb{C}}
\def\bbZ{\mathbb{Z}}
\def\bbN{\mathbb{N}}

\def\1{\mathbf{1}}

\def\bbC{\mathbb{C}}
\def\bbZ{\mathbb{Z}}
\def\bbN{\mathbb{N}}

\newcommand{\de}{\bigtriangleup}

\theoremstyle{plain}
\newtheorem{thm}{Theorem}[section]
\newtheorem{cor}[thm]{Corollary}
\newtheorem{lem}[thm]{Lemma}

\theoremstyle{definition}

\def\dss{\displaystyle}
\theoremstyle{remark}

\numberwithin{equation}{section}

\begin{document}

\def\dss{\displaystyle}

\title{Block Toeplitz Matrices: Multiplicative Properties}
\author{Muhammad Ahsan Khan}
\address{ Department of Mathematics, University of Kotli Azad Jammu $\&$ Kashmir, Kotli 11100, Azad Jammu $\&$ Kashmir, Pakistan}
\email{ahsankhan388@hotmail.com}
\begin{abstract}
 Given $A,B,C,$ and $D,$ block Toeplitz matrices, we will prove the necessary and sufficient condition for $AB-CD=0$, and $AB-CD$ to be a block Toeplitz matrix. In addition, with respect to change of basis, the characterization of normal block Toeplitz matrices with entries from  the algebra of diagonal matrices  is also obtained. 	
\end{abstract}
\maketitle
\section{Introduction}
A scalar Toeplitz matrix is an $n\times n$ matrix with the following structure: 
\[
A=
\begin{pmatrix} 
a_{0} &a_{1}&  a_{2}&\ldots  & a_{n-1}\\
a_{-1}&a_{0}       &  a_{1} &\ldots  & a_{n-2}\\
a_{-2}& a_{-1}     &  a_{0}        &\ldots  & a_{n-3}\\
\vdots    & \vdots      & \vdots        & \ddots\\
a_{1-n}&  a_{2-n} &  a_{3-n}   & \ldots & a_{0}
\end{pmatrix}
.\]
The entries depend upon the difference $i-j$ and hence they are constant down all the diagonals. The subject is many decades old; among monograph dedicated to the subject are \cite{grenand-szego,ioh} and \cite{widom}. These matrices are ubiquitous and arise naturally in several theoretical and applicative ﬁelds.  In particular, mathematical modelling of all the problems where to some extent the shift invariance appears in terms of space or of time. This shift invariance is contemplated in the structure of the matrix itself where a shift on south-eastern leaves the matrix fixed. For two dimensional problems the structure of Toeplitz may also occur block wise.

There is an immense literature concerning scalar Toeplitz matrices; on the other side, research related to block Toeplitz matrices is rather sparse. In this paper, we intend to add some remarks and results to the latter. Our starting point is the paper of Gu and Patton \cite{gu-patton}, which provides various properties of usual scalar Toeplitz and Hankel matrices; most of the results therein refer to products of these structured matrices, which, in general, are no more structured.

The area of block Toeplitz matrices is less studied, one of the reason being the new difficulties that appear with respect to the scalar case. Besides its theoretical interests, the subject is also important to multivariate control theory. We refer the reader to \cite{botsil,GGK,MAK,MAKDT,shalom} and further
references thereafter for a detailed study about block Toeplitz matrices. 

In \cite{gu-patton}, the authors have proved a variety of algebraic results about scalar Toeplitz matrices. Given Toeplitz matrices $A,B,C$, and $D$, their main result determines if the matrix $AB-CD$ is Toeplitz. The necessary and sufficient condition is a rank two matrix equation involving tensor products of the vectors defining $A,B,C,$ and $D$. They have also proved the necessary and sufficient condition for $AB-CD=0$. In addition to that, they have also completely characterized  normal Toeplitz matrices. The characterization of normal Toeplitz matrices has been discussed in \cite{FKKL,DIC,I}.  

The purpose of the present paper is to generalize some of the results of \cite{gu-patton} concerning the product of block Toeplitz matrices by way of introducing the special structure of the displacement matrix of a block Toeplitz matrix.

The outline of the paper is as follows: Notations and some basic facts about displacement matrices are presented in section~2. Then as we will be interested in block Toeplitz matrices, some basic properties concerning the product of these matrices are derived in section~3. Section ~4 is devoted to the study of normal block Toeplitz matrices with entries from the algebra of diagonal matrices. The last section is concerned with the applications of block Toeplitz matrices. 
\section{Preliminaries}
 As usual $\bbN$, $\bbC,$ and $\bbZ$ denote the set of natural numbers, complex numbers and integers respectively.  We symbolize by $\MM_{n}$ the algebra of $n\times n$ matrices and by $\DD_d$ the algebra of $d\times d$ diagonal matrices with entries from $\bbC$ . Throughout in this paper, we label
the indices from $0$ to $n-1$; so $A\in\MM_{n}$ is written $A=(a_{i,j})_{i,j=0}^{n-1}$ with $a_{i,j}\in\bbC$. Then $\TT_{n}\subset\MM_{n}$ is the space of scalar Toeplitz matrices $A = (a_{i-j})_{i,j=0}^{n-1}$. 
We will mostly be interested in block matrices, i.e., matrices whose elements are
not necessarily scalars, but elements in $\MM_{d}$. Thus a block Toeplitz matrix is actually an $ nd\times nd $ matrix, but which has been decomposed in $ n^2 $ blocks of dimension $ d $, and  these blocks are constant parallel to the main diagonal. In the sequel, we will use the following notations:
\begin{itemize}
	\item $ \MM_{n}\otimes \MM_d $ is the collection of $ n\times n $ block matrices whose entries all belong to $ \MM_d $;
	\item $ \TT_{n}\otimes\MM_d $ is the collection of $ n\times n $ block Toeplitz matrices whose entries all belong to $ \MM_d $;
	
	\item $ \DD_{n}\otimes\MM_{d}$ is the collection of $ n\times n $ diagonal block Toeplitz matrices whose entries all belong to $ \MM_d $; 
	
	\item $\CC_1\otimes\MM_{d}$ is the collection of all $n\times 1$ block matrices whose entries all belong to $\MM_d$;
	\item $\RR_1\otimes\MM_{d}$ is the collection of all $1\times n$ block matrices whose entries all belong to $\MM_d$. 
	
\end{itemize}
Obviously $\DD_{n}\otimes\MM_d\subset \TT_{n}\otimes\MM_d\subset \MM_{n}\otimes \MM_d $. For block diagonal matrices, we will use the notation 
\[
diag
\begin{pmatrix}
A_{1}&A_{2}&\cdots& A_{n}
\end{pmatrix}=
\begin{pmatrix}
A_{1} &  0     &\ldots  & 0\\
0     &  A_{2} &\ldots  & 0\\
\vdots& \vdots & \ddots & \vdots\\
0&  0  & \ldots & A_n
\end{pmatrix}
.\]

In most cases, it will suffice to consider block Toeplitz matrices with zero diagonals. In other cases $\tilde{A}+\tilde{A_0}$ will describe the most general block Toeplitz matrix , where $\tilde{A}$ is a block Toeplitz matrix with $0$ on the main diagonal and $\tilde{A_{0}}$ is the diagonal block  Toeplitz matrix.\\
If 
$a=
\begin{pmatrix}
	0\\
	a_{-1}\\
	\vdots\\
	a_{1-n}
\end{pmatrix}
$, then we define
$\hat{a}=
\begin{pmatrix}
	0\\
	{a_{1-n}}\\
	\vdots \\
	{a_{-1}}
\end{pmatrix}
$
and if  $b
=\begin{pmatrix}
	0\\
	a_1\\
	\vdots \\
	a_{n-1}
\end{pmatrix}$, then $\overline{\hat{b}}=
\begin{pmatrix} 
	0\\
	\overline{a_{n-1}}\\
	\vdots\\
	\overline{a_{1}}
\end{pmatrix}
.$ 
Let $I\in\MM_{d}$, be the identity matrix and  $S\in\MM_{n}\otimes\MM_d$  consisting of zero matrices, except for $I's$ along the subdiagonal, i.e., 
 \[
 S=
 \begin{pmatrix} 
 	0 &0&  0&\ldots  & 0&0\\
    I&0     & 0 &\ldots  & 0&0\\
 	0& I     &  0      &\ldots  & 0&0\\
 	\vdots    & \vdots      & \vdots    & \ddots&\vdots&\vdots\\
 	0&  0 &  0   & \ldots & I&0
 \end{pmatrix}
 .\]
 Then its adjoint is the matrix given by 
 \[
S^*=
\begin{pmatrix} 
	0 &I&  0&\ldots  & 0&0\\
	0&0     & I &\ldots  & 0&0\\
	0& 0     &  0      &\ldots  & 0&0\\
	\vdots    & \vdots      & \vdots    & \ddots&\vdots&\vdots\\
		0&  0 &  0   & \ldots &0 &I\\
	0&  0 &  0   & \ldots & 0&0
\end{pmatrix}
.\] 
Note that $S^n=S^{*n}=0.$

 If we view $S$ as a linear operator acting on the space $\CC_{1}\otimes \MM_{d}$. Then $S$ shifts the components of a column vector one position down, with a zero matrix appearing in the first position. While its adjoint $S^*$ shifts the components of a column vector one position up, with a zero matrix appearing in the last position. 
For any $M\in\MM_{n}\otimes\MM_{d}$, the displacement matrix is defined as
 \[
 \de(M):=M-SMS^{\ast}
 .\]
We will use this matrix to determine whether the difference of the matrix products is Toeplitz.  For other kinds of displacement matrices see \cite{HK} and \cite{KS} . We denote the vector
 $\begin{pmatrix}
 	I\\
 	0\\
 	\vdots\\
 	0
 \end{pmatrix}^T\in\RR_1\otimes\MM_{d}$ by $P_+$, then its adjoint is the vector $P_+^\ast=\begin{pmatrix}
 I\\
 0\\
 \vdots\\
 0
\end{pmatrix}\in\CC_1\otimes\MM_{d}$.

The following lemma is quite useful for proving our main results of section~ 3.
\begin{lem}\label{dis}
	If $M\in\MM_{n}\otimes\MM_{d}$, then  $M=\dss\sum_{k=0}^{n-1}S^k(\de(M))S^{*k}$.
\end{lem}
\begin{proof}
\begin{align*}
\dss\sum_{k=0}^{n-1}S^k(\bigtriangleup(M))S^{*k}
&=\dss\sum_{k=0}^{n-1}S^k(M-SMS^*)S^{*k}\\
&=\dss\sum_{k=0}^{n-1}(S^kMS^{*k}-S^{k+1}MS^{*k+1})\\
&=M-S^nMS^{*n}=M
\end{align*}
\end{proof}
Thus if want to show that $M=0$, then it will  sufficient to show that  $\de(M)=0$\\
.

\section{Block Toeplitz product}
 In this section, we will generalized some important results of \cite{gu-patton}. 
The following Lemma describes the necessary and a sufficient condition for a block matrix $A$ to be a block  Toeplitz matrix. 
\begin{lem}\label{Toeplitz}
	$A\in\MM_{n}\otimes\MM_{d}$  is Toeplitz if and only if there exist  $X,X^{\prime}\in\CC_1\otimes\MM_{d}$  such that, $\bigtriangleup(A)=XP_{+}+P_{+}^{\ast}X^{\prime\ast}.$
\end{lem}
\begin{proof}
	Suppose that $A=(A_{i-j})_{i,j=0}^{n-1}\in\TT_{n}\otimes\MM_{d}$. Since, the displacement matrix for $A$ is defined as $\de(A)=A-SAS^{\ast}.$
	Then  simple computation yields that 
	\[\de(A)=
	\begin{pmatrix} 
	A_{0} & A_{1} & A_{2}&\ldots  & A_{n-1}\\
	A_{-1}&  0    &  0   &\ldots  & 0\\
	A_{-2}&  0    &  0   &\ldots  & 0\\
	\vdots& \vdots   & \vdots & \ddots\\
	A_{-(n-1)}&  0 &  0  & \ldots & 0
	\end{pmatrix}
	\] 
	If we take 
	$
	X = \begin{pmatrix}
	A_{0}\\
	A_{-1}\\
	\vdots\\
	A_{1-n}
	\end{pmatrix}
	$ and $X^{\prime} =
	\begin{pmatrix}
	0 \\
	A_{1}\\
	\vdots\\
	A_{n-1}
	\end{pmatrix}
	$, then one can easily check that $$\de(A)=XP_{+}+P_{+}^{\ast}X^{\prime\ast}.$$\\
	For the converse, let $A=(A_{i,j})_{i,j=0}^{n-1}\in\MM_{n}\otimes\MM_{d}$. 
	Suppose then  that $
X = \begin{pmatrix}
	X_{0}\\
	 X_{1}\\
	 \vdots\\
	  X_{n-1}
	\end{pmatrix}
	$
	and $X^{\prime}=
	\begin{pmatrix}
	X_{0}^{\prime}\\
	X_{1}^{\prime}\\
	\vdots\\
	 X_{n-1}^{\prime}
	 \end{pmatrix}
	 $ be  vectors in $\CC_1\otimes\MM_{d}$, since we have 
	\[
	\de(A)=XP_{+}+P_{+}^{\ast}X^{\prime\ast}
	\] $\implies$  
	\[
	A=SAS^{\ast}+XP_{+}+P_{+}^{\ast}X^{\ast\prime}
	\]$ \implies$ 
	\[A=
	\begin{pmatrix} 
	X_{0}+X_{0}^{\prime\ast}&X_{1}^{\prime\ast}&X_{2}^{\prime\ast}&\ldots&X_{n-1}^{\prime\ast}\\
	X_{1} &  A_{0,0} &  A_{0,1} &\ldots & A_{0,n-2}\\
	\vdots & \vdots  & \vdots  & \ddots\\
	X_{n-1}&   A_{n-2,0} &  A_{n-1,1} & \ldots & A_{n-2,n-2}
	\end{pmatrix}
	\]  
	Compairing corresponding entries yields
	$A_{i_{1},j_{1}}=A_{i_{2},j_{2}}$, whenever $i_{1}-j_{1}=i_{2}-j_{2}$, where $0\leq i_{1}, i_{2},j_{1},j_{2}\leq n-1$., i.e., $A$ is a block Toeplitz matrix.  
	
\end{proof}
The space $\TT_{n}\otimes \MM_{d}$ do not form an algebra , so it is interesting to study about the algebraic properties of $A`s$ in $\TT_{n}\otimes\MM_{d}$.

In the rest of this section, if $A=(A_{i-j})_{i,j=0}^{n-1}$ is any block Toeplitz matrix  then for simplification we write $A=\tilde{A}+\tilde{A_{0}}$, where $\tilde{A}\in\TT_{n}\otimes\MM_d$, with $0$ on the main diagonal and $\tilde{A_0}\in\DD_{n}\otimes\MM_{d}$.  
 The next lemma is the main technical result of this section. 
\begin{lem}\label{product}
	Let  $\tilde{A_{-}}=
	\begin{pmatrix}
		0\\
		 A_{-1}\\
		 \vdots\\
		  A_{-(n-1)}
	\end{pmatrix}
	$,
	$
	\tilde{B_{+}}=
	\begin{pmatrix}
		0\\
		 B_{1}\\
		 \vdots \\
		 B_{n-1}
	\end{pmatrix}^T
	$,
	$
	\tilde{A_{+}}=
	\begin{pmatrix}
		0\\
		 A_{n-1}\\
		 \vdots\\
		  A_{1}
	\end{pmatrix}
	$
	and 
	$
	\tilde{B_{-}}=
	\begin{pmatrix}
		0\\
		 B_{1-n}\\
		 \vdots\\
		  B_{-1}
	\end{pmatrix}^T
	$ be vectors in  $\CC_{1}\otimes\MM_{d}$. Suppose $A=\tilde{A}+\tilde{A_{0}}$ and $B=\tilde{B}+\tilde{B_{0}}$, then there exist vectors $Y$, and $Y^\prime$ in $\CC_1\otimes\MM_{d}$, such that
	$
	\bigtriangleup(AB)=YP_{+}+P_{+}^{\ast}Y^{\prime\ast}+\tilde{A}_{-}\tilde{B}_{+}-\tilde{A}_{+}\tilde{B}_{-}.
	$

\end{lem}
\begin{proof}
	We have
	\begin{align*}
	\bigtriangleup(AB)
	&=\bigtriangleup[(\tilde{A}+\tilde{A_{0}})(\tilde{B}+\tilde{B_{0}})]\\
	&=\bigtriangleup[\tilde{A}\tilde{B}+\tilde{A}\tilde{B_{0}}+\tilde{A_{0}}\tilde{B}+\tilde{A_{0}}\tilde{B_{0}}]\\
	&=\bigtriangleup(\tilde{A}\tilde{B})+\de(\tilde{A}\tilde{B_{0}})+\de(\tilde{A_{0}}\tilde{B})+\de(\tilde{A_{0}}\tilde{B_{0}})  
	\end{align*}
	Since $\tilde{B_{0}}\in\DD_{n}\otimes\MM_{d}$, then $\tilde{A}\tilde{B_{0}}\in\TT_{n}\otimes\MM_d$, therefore by Lemma \ref{Toeplitz}, there exist $U, U^\prime\in\CC_1\otimes\MM_{d}$ such that $\de(\tilde{A}\tilde{B_{0}})=UP_++P_+^*U^{\prime\ast}$. Similarly $\de(\tilde{A_{0}}\tilde{B})= VP_{+}+P_{+}^{\ast}V^{\prime\ast}$, and $\de({\tilde{A_{0}}\tilde{B_{0}}})=WP_{+}+P_{+}^{\ast}W^{\prime\ast}$, with  $V,V^\prime,W, W^{\prime}\in\CC_1\otimes\MM_{d} $.
	Therefore we can write 
	\begin{equation}\label{del1}
	\bigtriangleup(AB)=\de(\tilde{A}\tilde{B})+(U+V+W)P_{+}+(U^{\prime\ast}+V^{\prime\ast}+W^{\prime\ast})P_{+}^{\ast}.
	\end{equation}
	We have entry at the position $(i,j)$ of $\tilde{A}\tilde{B}$ is 
	\begin{equation}\label{product formula}
	(\tilde{A}\tilde{B})_{i,j}=\sum_{k=0}^{n-1}A_{k-i}B_{j-k}, \quad 0\leq i,j\leq n-1.
	\end{equation} 
	Using the formula \eqref{product formula} and definition of $S$, one obtains
	\begin{equation}\label{eq1}
(\tilde{A}\tilde{B}-S\tilde{A}\tilde{B}S^{\ast})_{i,j}=	\begin{cases}
	A_{-i}B_{j}-A_{n-i}B_{j-n} \quad\hbox{for}\quad  1\leq i,j\leq n-1,\\
	(\tilde{A}\tilde{B})_{0,j}\quad \hbox{for}\quad 1\leq j\leq n-1,\\
	(\tilde{A}\tilde{B})_{i,0}\quad  \hbox{for}\quad 0\leq i\leq  n-1.
	\end{cases}
	\end{equation}
	Therefore 
	\begin{equation}\label{del2}
	\de(\tilde{A}\tilde{B})=XP_{0}+P_{0}^{\ast}X^{\prime\ast}+\tilde{A}_{-}\tilde{B}_{+}-\tilde{A}_{+}\tilde{B}_{-},
	\end{equation}
	where 
	$
	X=
	\begin{pmatrix}
	(\tilde{A}\tilde{B})_{0,0}\\
	(\tilde{A}\tilde{B})_{1,0}\\
	\vdots \\
	(\tilde{A}\tilde{B})_{n-1,0}
	\end{pmatrix}
	$,
	$
	X^{\prime}=
	\begin{pmatrix}
	0\\
	(\tilde{A}
	\tilde{B})_{0,0}\\
	\vdots\\
	(\tilde{A}
	\tilde{B})_{0,n-1}
	\end{pmatrix}
	 .$ 
	Combining \eqref{del1} and \eqref{del2} yields
	\begin{align*}
	\bigtriangleup(AB)
	&=(U+V+W+X)P_{+}+(U^{\prime\ast}+V^{\prime\ast}+W^{\prime\ast}+X^{\prime\ast})P_{+}^{\ast}+\tilde{A}_{-}\tilde{B}_{+}-\tilde{A}_{+}\tilde{B}_{-}\\
	&=YP_{+}+P_{+}^{\ast}Y^{\prime\ast}+\tilde{A}_{-}\tilde{B}_{+}-\tilde{A}_{+}\tilde{B}_{-},
	\end{align*}
	where $Y=U+V+W+X_{0}$ and $Y^{\prime\ast}=U^{\prime\ast}+V^{\prime\ast}+W^{\prime\ast}+X_{0}^{\prime\ast}$.
\end{proof}
The next Theorem is the most important result of this paper. It gives answer to the question that when the product $AB-CD$ of  block Toeplitz matrices $A,B,C,$ and $D$ is in $\TT_{n}\otimes\MM_{d}.$ 
\begin{thm}\label{Main}
	 Let $\tilde{A}_-, \tilde{A_{+}},\tilde{C_{-}},\tilde{C_{+}}\in\CC_1\otimes\MM_{d}$ and $ \tilde{B_{+}},\tilde{B_{-}},\tilde{D_{+}},\tilde{D_{-}}\in\RR_1\otimes\MM_{d}$ with $0$ in the zeroth component.
	 Suppose that $A=\tilde{A}+\tilde{A_{0}}$, $B=\tilde{B}+\tilde{B_{0}}$, $C=\tilde{C}+\tilde{C_{0}}$ and $D=\tilde{D}+\tilde{D_{0}}$, then $AB-CD\in\TT_{n}\otimes\MM_{d}$ if and only if  $\tilde{A_{-}}\tilde{B_{+}}-\tilde{A_{+}}\tilde{B_{-}}=\tilde{C_{-}}\tilde{D_{+}}-\tilde{C_{+}}\tilde{D_{-}}$.\\

\end{thm}
\begin{proof}
	 By Lemma \ref{product}  
	\begin{align*}
	\de(AB-CD)
	&=\de(AB)-\de(CD)\\
	&=(Y-Z)P_++P_+^\ast (Y^{\prime\ast}-Z^{\prime\ast})+ \tilde{A_{-}}\tilde{B_{+}}-\tilde{A_{+}}\tilde{B_{-}}-\tilde{C_{-}}\tilde{D_{+}}+\tilde{C_{+}}\tilde{D_{-}}
	.\end{align*}
One can check that  the last four terms on the right side of the above equation involve vectors with $0$ in the zeroth component. Lemma \ref{Toeplitz} imply that, $AB-CD$ is in $\TT_{n}\otimes\MM_{d}$ if and only if $\tilde{A_{-}}\tilde{B_{+}}-\tilde{A_{+}}\tilde{B_{-}}=\tilde{C_{-}}\tilde{D_{+}}-\tilde{C_{+}}\tilde{D_{-}}$.\\
	
\end{proof}
\begin{cor}
Let $\tilde{A_{-}}, \tilde{A_{+}}\in\CC_1\otimes\MM_{d}$, and $\tilde{B_{+}},\tilde{B_{-}}\in\RR_1\otimes\MM_{d}$ with $0$	in the zeroth component. Let $A=\tilde{A}+\tilde{A_0}$ and $B=\tilde{B}+\tilde{B_{0}}$, then $AB\in\TT_{n}\otimes\MM_{d}$ if and only if $\tilde{A_{-}}\tilde{B_{+}}-\tilde{A_{+}}\tilde{B}_{-}=0$.
\end{cor}
\begin{proof}
The proof follows immediately from Theorem \ref{Main} by taking $C=D=0$.
\end{proof}
\begin{thm}
	 Let $Y,Y^\prime$, $Z$ and $Z^\prime$ be 
	vectors in $\CC_1\otimes\MM_{d}$ with $0$ on the zero component. Let $AB-CD\in\TT_{n}\otimes\MM_{d}$, then $AB = CD$ if and only if $Y=Z$ and $Y^\prime=Z^\prime.$
\end{thm}
\begin{proof}
	Let $AB-CD\in\TT_{n}\otimes\MM_{d}$, then by Lemma \ref{Toeplitz}, there exist $Y,Y^\prime$ $Z$ and $Z^\prime$ in $\CC_1\otimes\MM_{d}$, such that $$\de(AB-CD)=(Y-Z)P_++P_+^\ast (Y^{\prime\ast}-Z^{\prime\ast}).$$
	It follows from Lemma \ref{dis}, that $AB=CD$ if and only if $\de(AB-CD)=0$.  The latter equation holds if and only if the vectors which form the product with $P_+$ are $0$. That is, $AB=CD$ if and only if $Y=Z$ and $Y^\prime=Z^\prime$.
\end{proof}
\section{Normal block Toeplitz matrices}
Normal matrices are matrices that include unitary matrices and enjoy several of the same properties as unitary matrices. The general problem of characterizing normal block Toeplitz matrices with entries from $\MM_{d}$ is an open problem. As a first step, we take entries of block Toeplitz matrices from $\DD_{d}$, then  we will obtain the characterization of normal  block Toeplitz matrices. We start with the following result, which is Lemma 5.1 of \cite{MAK} .
\begin{lem}\label{Change of basis}
	Suppose $A=(A_{i-j})_{i,j=0}^{n-1}\in\TT_{n}\otimes\DD_{d}$. Then there is a change of basis that brings $A$ into the following form
	\[A^{\prime}=diag
	\begin{pmatrix}
		A_{1}^{\prime}&A_{2}^{\prime}&\cdots A_{d}^{\prime}
	\end{pmatrix}
	,\]
	where for every $k=1,2,\cdots d$, $A_{k}^{\prime}\in\TT_{d}$.
\end{lem}
The next result from \cite{gu-patton} characterized  normal  Toeplitz matrices among all Toeplitz matrices.  

\begin{thm}\label{scalar normal}
	Let $A=(a_{i-j})_{i,j=0}^{n-1}\in\TT_{n}$, then $A$ is normal if and only if  either $a=\lambda \hat{b}$  for some $|\lambda|=1$ or $a=\lambda\bar{\hat{b}}$ for some $|\lambda|=1$.
\end{thm} 
The following result is the main result of this section. Although, it is a trivial implications of Lemma \ref{Change of basis} but it will provide a basis for future research concerning the  characterization of normal block Toeplitz matrices with entries from $\MM_{d}$.
\begin{thm}
	Let $A\in\TT_n\otimes\DD_{d}$, then $A$ is normal if and only if there exist scalars $\lambda_{1},\lambda_{2},\cdots\lambda_{d}$ with $|\lambda_{k}|=1$, such that for every $k=1,2,\cdots d$, either $a_{k}=\lambda_{k}\hat{b} $ or $a_{k}=\lambda_k\overline{\hat{b}}_{k}$. 
\end{thm}
\begin{proof}
	Suppose $A\in\TT_{n}\otimes\DD_{d}$, then by Lemma \ref{Change of basis}, $A$ has the form
	\[A^{\prime}=diag
	\begin{pmatrix}
		A_{1}^{\prime}&A_{2}^{\prime}&\cdots A_{d}^{\prime}
	\end{pmatrix}
	,\]
	where for every $k=1,2,\cdots d$, $A_{k}^{\prime}=(a_{r-s,k})_{r,s=0}^{d}$.  
	Since $A^\prime\in\DD_{n}\otimes\MM_{d}$, then  $AA^*=A^*A$ if and only if $A_{k}^\prime A_{k}^{\prime *}=A_{k}^{\prime *} A_{k}^\prime$. By Theorem \ref{scalar normal}, each $A_k^\prime$ is normal if and only if there exist scalars   $\lambda_{1},\lambda_{2},\cdots\lambda_{d}$, with $|\lambda_{k}|=1$, such that  either $a_k=\lambda_k \hat{b}$  or $a_k=\lambda_k\overline{\hat{b}}_{k}$, where 
	$a_{k}=
	\begin{pmatrix}
		a_{-1,k}\\
		a_{-2,k}\\
		\vdots \\
		a_{-d,k}
	\end{pmatrix}
$, and $b_{k}=
	\begin{pmatrix}
		a_{1,k}\\
		a_{2,k} \\
		\vdots \\
		a_{d,k}
	\end{pmatrix}
	$, $k=1,2,\cdots d$.
The proof is therefore complete.
\end{proof}

\section{Applications of block Toeplitz matrices}
Block Toeplitz matrices appear in various theoretical and applicative fields. By means of this section, we want to make sure readers the importance of the study of these matrices.  
 \subsection{Problems modelled by (block) Toeplitz matrices}
 There is a natural relation between problems involving Toeplitz structures and power series  (even with Laurent series in case of infinite Toeplitz structure) which allows one to shift from algorithms for Toeplitz computations to algorithms for power series computations and vice versa.  Fast Fourier Transform in this way becomes a principal tool for all the computations involving Toeplitz-like matrices.
 In the last decades a large amount of research has been concentrated on the analysis of algorithms for Toeplitz matrices. In particular, there are two specific problems:
 \begin{itemize}
 	\item[(i)] solving  linear system involving  Toeplitz matrix of dimension $d$, where the matrix is generated by a real function; see \cite{B, F,GR}.
 	\item [(ii)] finding vector invariant under the action of an infinite block Toeplitz matrix; see \cite{B,GR, RD}.
 \end{itemize}
 
 \subsection{Problem in Queueing Theory}
 An infinite block Toeplitz matrix in Hessenberg form, muddled more exquisite problems. This is because of an infinite features of the problem. Currently there is no direct method exists, and the most appropriate solution techniques are dependent on finding the fixed points of the matrix (by viewing matrix as an operator). Such type of  methods utterly use the Toeplitz structure without utilizing the acceleration permitted by the Fast Fourier Transform (for details see \cite{B}).
 \subsection{Image restoration}
 There is another interesting appllication of block Toeplitz matrices related to  blurring and deblurring models in digital image restoration. Suppose that the blur of a single point of an image does not depend on  the location of the point, that is, it is shift invariant, and is deﬁned by the Point-Spread Function (PSF). Such type of a  function has compact support, indeed, a point is blurred into a small fleck of light with dark everywhere except that in a small neighborhood of the point.
 The relation between the blurred and noisy image, stored as a vector $B$ and the real image, represented by a vector $X$ has the form
 \[
 AX=B-noise.
 \]
 Because of the shift invariance of the PSF, $A$ is a block Toeplitz matrix  with Toeplitz blocks. Due to the local eﬀect of the blur, the PSF has compact support so that $A$ is block banded with banded blocks. Conventionally, $A$ is ill conditioned so that by solving the system $AX=B$ we get by ignoring the noise provides a highly perturbed solution. For example the PSF which transforms a unit point of light into 
 the $3\times 3$
 square
 \[
 \frac{1}{15}
 \begin{pmatrix}
 	1& 2& 1\\
 	2& 3& 2\\
 	1& 2& 1
 \end{pmatrix}\]
 
 yields the following block Toeplitz matrix
 
 \[
 A=\frac{1}{15}
 \begin{pmatrix}
 	U& T& & \\
 	
 	T&U&T& \\
 	& \ddots&\ddots&\ddots\\
 	&&T& U&T\\
 	&&&T&U	
 \end{pmatrix}
 ,\]
 where $T=\begin{pmatrix}
 	2& 1& & \\
 	
 	1&2&1& \\
 	& \ddots&\ddots&\ddots\\
 	&&1& 2&1\\
 	&&&1&2	
 \end{pmatrix}$ and $U=\begin{pmatrix}
 	3& 2& & \\
 	
 	2&3&2& \\
 	& \ddots&\ddots&\ddots\\
 	&&2& 3&3\\
 	&&&2&3	
 \end{pmatrix}.$ 
 Restoring a blurred image in this way, is reduced to solving a block banded block Toeplitz systems with banded Toeplitz blocks. Such type of basic problems appears in many forms in image processing \cite{ HNO,RHE, NCT, WNC}.
 \subsection{Block Toeplitz matrices in signal processing}
 Signal processing plays a vital role in multiple industries and it is the technology of future. We now give two interesting examples of block Toeplitz matrices
 that are frequently used in Signal Processing as well as in  Communications and  Information Theory; see \cite{GR,JV, WK}. \\
 (i) Let $\{\eta_d\}_{d\in\bbZ}$ be a sequence in $\CC_1\otimes\bbC$
 , i.e., $\eta_d\in\CC_1\otimes\bbC$ for every $d\in\bbZ$. If $j,k\in\bbZ$ with $j \leq k$, then
 \[
 \eta_{k:j}=\begin{pmatrix}
 	\eta_k\\
 	\eta_{k-1}\\
 	\vdots\\
 	\eta_j
 \end{pmatrix}
 \]
 Then consider a discrete time causal ﬁnite impulse response (FIR) multiple input multiple output (MIMO) ﬁlter, that is, a ﬁlter given as
 \begin{equation}\label{fil}
 	\zeta_d=\sum_{k=0}^{r}A_{-k}\eta_{d-k}\quad \hbox{for every}\quad d\in\bbZ,
 \end{equation}
 where the filters taps $A_{-k}$, with $0\leq k\leq r$ are $m\times n$ block matrices, and the input $\{\eta_d\}_{d\in\bbZ}$ and the output $\{\zeta_d\}_{d\in\bbZ}$ of the filters satisfy that $\eta_d\in\CC_1\otimes\bbC$ and $\zeta_d\in\CC_1\otimes \bbC$ for every $d\in\bbZ$. It follows then from  (\ref{fil}) that 
 \[
 \zeta_{d:1}=A_d \eta_{d:1-r}, \quad \hbox{for every } d\in\bbZ,
 \]	where $A_d$ is the $d \times (d + r)$ block Toeplitz matrix with $m\times n$ blocks given by $A_d=(A_{j-k})_{1\leq j\leq d, 1\leq k\leq d+r}$ with $A_{j-k}=0_{m\times n}$ when $j-k\neq -d,\cdots -1,0$. Thus, representations of discrete time causal FIR MIMO ﬁlters in terms of matrix are block Toeplitz.
 
 (ii) We now consider a vector wide sense stationary (WSS) process. Suppose that $\eta_n$ be a random vector of dimension $d$ for every $n\in\bbZ$. Suppose that the $d$ dimensional multivariate random processes $\{\eta_n\}_{n\in\bbZ}$  has constant mean, i.e., 
 \[
 E(\eta_{n_1})=E(\eta_{n2})\quad \hbox{for every }\quad n_1,n_2\in\bbZ. 
 \]
 and if it obeys
 \begin{equation}\label{random}
 	E(\eta_{n_1}\eta_{n_2}^*)=E(\eta_{n_3}\eta_{n_4}^*)
 \end{equation}
 whenever $n_1-n_2=n_3-n_4.$ $(\ref{random})$ implies that
 \begin{equation*}
 	E(\eta_{n:1}\eta_{n:1}^*)=(A_{j-k})_{j,k=1}^{n} \quad\hbox{for any} \quad n\in\bbN,
 \end{equation*}
 where $A_{j-k}=E(\eta_k\eta_j^*)\in\TT_{n}\otimes\MM_{d}$, for every $j,k\in\bbN$.
 Thus, the corelation matrix $E(\eta_{n:1}\eta_{n:1}^*)$ of the random of vector $\eta_{n:1}$ is
 an $n\times n$ block Toeplitz matrix with $d \times d$ blocks for every $n\in\bbN$.
\section*{Acknowledgments}  The author would like to thank the referees for their comments, which leads to an improved version of the paper.

\end{document}